\documentclass[10pt]{amsart}
\usepackage{amsmath}
\usepackage[all]{xy}
\usepackage{amssymb} 
\usepackage{epsbox}
\usepackage{amsthm}
\usepackage{graphicx}

\begin{document}  
\title{A functor-valued extension of knot quandles}
\author{Tetsuya Ito}
\address{Graduate School of Mathematical Science, University of Tokyo, 3-8-1 Komaba, Meguro-ku, Tokyo, 153-8914, Japan}
\email{tetitoh@ms.u-tokyo.ac.jp}
\subjclass[2010]{Primary~57M25, Secondary~57M27}
\keywords{quandle, knot, quandle homology, cocycle invariant, quandle invariant functor}

\def\hoop{\unitlength.1em
  \begin{minipage}{9\unitlength}
    \begin{picture}(30,20)
      \put(15,10){\circle{10}}
      \put(15,10){\circle*{0.5}}
      \put(0,10){\line(1,0){9,5}}
    \end{picture}
  \end{minipage}
  \;\;\;\;
}

 
\newtheorem{thm}{Theorem}
\newtheorem{cor}{Corollary}
\newtheorem{lem}{Lemma}
\newtheorem{prop}{Proposition}
 {\theoremstyle{definition}
  \newtheorem{defn}{Definition}}
 {\theoremstyle{definition}
  \newtheorem{exam}{Example}}
 {\theoremstyle{definition}
  \newtheorem{rem}{Remark}}
  {\theoremstyle{definition}
  \newtheorem{ques*}{Question}}
 \newcommand{\os}{ \: \overline{*} \:}
 \newcommand{\Z}{ \mathbb{Z}}
 \newcommand{\Q}{ \mathbb{Q}}
 \newcommand{\R}{ \mathbb{R}}
 \newcommand{\mA}{ \mathcal{A}}
  \newcommand{\mSA}{ \mathcal{SA}}
  \newcommand{\aut}{ \textrm{Aut}\,}
  \newcommand{\Inn}{ \textrm{Inn}\,}
 \newcommand{\Ass}{\textrm{Ass}\,}
  \newcommand{\Hom}{\textrm{Hom}}
 

\begin{abstract}
For an oriented knot $K$, we construct a functor from the category of pointed quandles to the category of quandles in three different ways. We also extend the quandle cocycle invariants of knots by using these quandle-valued invariant of knots, and study their properties. 
\end{abstract}
\maketitle

\section{Introduction}

 A quandle is a set $Q$ with a binary operation $*$ which satisfies some axioms, and the pointed quandle is a pair $(Q,h)$ consisting of a quandle $Q$ and its element $h$.
 For an oriented knot $K$, one can associate a quandle $Q_{K}$ called the {\it knot quandle}. The knot quandle distinguishes all knots up to orientation. Moreover, using the homology theories of quandles, the knot quandle provides a knot invariant called a {\it quandle cocycle invariant} \cite{cegs}, \cite{cjkls}.

The aim of this paper is to extend the knot quandle as a functor.
We construct the {\it quandle invariant functor} $I_{K}$ for each oriented knot $K$, which is a functor from $\mathcal{PQ}$, the category of pointed quandles, to $\mathcal{Q}$, the category of quandles. Thus, for each pointed quandle $(Q,h)$ we obtain a quandle-valued invariant of a knot $I_{K}(Q,h)$. The classical knot quandle $Q_{K}$ appears as $I_{K}(T_{1})$, the quandle valued invariant which corresponds to the trivial $1$-quandle $T_{1}$. We also construct an extension of quandle cocycle invariants by using the quandle-valued invariant of knots.

Our results are summarized as follows. 
\begin{thm}
\label{thm:main}
Let $K$ be an oriented knot. Then there exists a functor $I_{K}: \mathcal{PQ} \rightarrow \mathcal{Q}$ from the category of pointed quandles to the category of quandles, having the following properties.
\begin{enumerate}
\item For the trivial $1$-quandle $T_{1}$, $I_{K}(T_{1})$ is the knot quandle $Q_{K}$.
\item $H_{1}^{Q}(I_{K}(Q,h);\Z) \cong H_{1}^{Q}(Q;\Z)$.
\item If $(Q,h)$ is a finite pointed quandle, then there exists a characteristic homology class $[K]_{Q,h} \in H_{2}^{Q}(I_{K}(Q,h); \Z)$ which vanishes if and only if $K$ is unknot.   
\item For the dual knot $-K^{*}$, which is the mirror of $K$ with the opposite orientation, $I_{K} =I_{-K^{*}}$ and $[-K^{*}]_{Q,h} = -[K]_{Q,h}$ hold.
\end{enumerate}
\end{thm}

Using the characteristic class $[K]_{Q,h}$, we define a cocycle invariant as in the classical knot quandles. The explicit construction of cocycle invariants will be given in section 6. In this section we also study the fundamental properties of this extension of cocycle invariants.

   We construct the quandle invariant functor by three different ways. 
The first method is algebraic. We use a representation of the braid groups derived from a pointed quandle, which can be seen as an generalization of the Artin representation of braid groups introduced in \cite{cp}.
  The second method is combinatorial and uses knot diagrams. We define the quandle invariants by giving a presentation, which extends the Wirtinger presentation of the knot quandles. This point of view is useful when we extend the quandle cocycle invariants.
  The last method is geometric. We construct a topological pair of space whose positive fundamental quandle coincides with the quandle invariant of knots. Such a spatial realization of the quandle invariant is obtained by gluing a topological space along the knot. 
    
 Our quandle-valued invariants are generalizations of the group-valued invariant defined by Crisp-Paris \cite{cp} and Wada \cite{w}. Indeed, many arguments in this paper is a direct generalization of the arguments in \cite{cp}. A new aspect which did not appear in \cite{cp} is the homology theory, which is used to define a cocycle invariant. 

Finally, we remark that although we restrict our attention to knots, but the constructions of the quandle invariant functors and the characteristic classes are valid for oriented links as well with some modifications. Thus our results are directly extended for link cases.\\

\textbf{Acknowledgment.}
 This research was supported by JSPS Research Fellowships for Young Scientists.
  
\section{Quandles, racks and braids}

In this paper we denote by $B_{n}$ the degree $n$ braid group, and by $\sigma_{1},\ldots,\sigma_{n-1}$ the standard generators of $B_{n}$. The closure of a braid $\beta$, which is an oriented link in $S^{3}$, is denoted by $\widehat{\beta}$.

A quandle is a set $Q$ with the binary operation $*$, which satisfies the following three axioms.
\begin{description}
\item[Q1] $a*a =a$ for all $a \in Q$.
\item[Q2] For any $a,b \in Q$ there exists a unique element $a \os b \in Q$ such that $a = (a\os b) * b = (a * b) \os b$.
\item[Q3] $(a*b)*c = (a*c)*(b*c)$ for all $a,b,c \in Q$. 
\end{description}

If $(Q,*)$ does not satisfy the axiom {\bf Q1} but satisfy both {\bf Q2} and {\bf Q3}, then the $(Q,*)$ is called a {\it rack}.

\begin{exam}
Here we present examples of quandle, which will be used in the paper.
\begin{itemize}
\item A set $X$ with the operation $*$ defined by $x*y = x$ is a quandle. We call this the {\it trivial quandle}. If the cardinal of $X$ is $n$, we call it the {\it trivial } $n$-{\it quandle} and denote it by $T_{n}$.
\item A group $G$ can be considered as a quandle by the operation $*$ defined by the conjugation $x*y = y^{-1}xy$. We call this quandle {\it the conjugacy quandle} associated to the group $G$ and denote by $Q_{G}$.
Conversely, for a quandle $Q$, one can obtain the {\it associated group} $\Ass(Q)$ defined by the presentation 
\[ \Ass(Q) = \langle q \in Q \:| \: p*q = q^{-1}pq \rangle. \] 
\end{itemize}
\end{exam}

We call a pair $(Q,h)$ consisting of a quandle $Q$ and its element $h \in Q$ a {\it pointed quandle}. 
A map between two quandles $\tau: (Q,*_{Q}) \rightarrow (P,*_{P})$ is called a {\it morphism} if $\tau$ preserves the operation $*$, that is, $\tau(a*_{Q} b) = \tau(a)*_{P}\tau(b)$ holds for all $a,b \in Q$. We denote by $\aut(Q)$ the group of automorphisms of $Q$. A morphism between pointed quandles $(Q,h)$ and $(P,i)$ is, by definition, a quandle morphism $f$ which satisfies $f(h)=i$. We denote the category of pointed quandles and the category of quandles by $\mathcal{PQ}$ and $\mathcal{Q}$ respectively.

As in the group case, the notion of the free quandles, free products, and the presentation of quandles are defined by the similar way. 

\section{Representations of the braid group associated to pointed quandles}

In this section we define a representation of the braid group $\rho_{Q,h}:B_{n} \rightarrow \aut(Q^{*n})$ for a pointed quandle $(Q,h)$.  
Let $Q^{*n} = Q_{1}*Q_{2}*\cdots *Q_{n}$ be the free product of $n$-copies of the quandle $Q$.
For an element $q \in Q$, we denote by $q_{i}$ the element in $Q_{i} \subset Q^{*n}$ which corresponds to $q$. 
For each integer $k=1,2,\ldots,n-1$, let $\tau_{k}$ be an automorphism of $Q^{*n}$ defined by
\[
\tau_{k}\: :
\left\{
\begin{array}{rll}
q_{k} & \mapsto & q_{k+1}\os h_{k} \\
q_{k+1} & \mapsto & q_{k} * h_{k} \\
q_{i} & \mapsto & q_{i} \;\;\;\; (i \neq k,k+1)
\end{array}
\right.
\]

The following proposition is easily confirmed by a direct calculation.
\begin{prop}
The map $\rho_{Q,h}: B_{n} \rightarrow \aut (Q)$ defined by $\rho_{Q,h}(\sigma_{i})= \tau_{i}$ is a group homomorphism.
\end{prop} 

We call this representation the {\it representation associated to a pointed quandle} $(Q,h)$.
By considering the associated group of $Q$, we also have a representation
$\rho'_{Q,h}: B_{n} \rightarrow \aut(\Ass(Q^{*n}))$, which is explicitly written as 
\[ 
\rho'_{Q,h}(\sigma_{i})\: :
\left\{
\begin{array}{rll}
q_{k} & \mapsto & h_{k}^{-1}q_{k+1}h_{k} \\
q_{k+1} & \mapsto & h_{k}q_{k}h_{k}^{-1} \\
q_{i} & \mapsto & q_{i} \;\;\;\; (i \neq k,k+1)
\end{array}
\right.
\]

 This representation is called the {\it Artin type representation} of $B_{n}$ associated to the pair $(\Ass(Q),h)$, which is defined in \cite{cp}.

\begin{exam}
\label{exam:classical}
Let $F_{n}$ be the rank $n$ free group generated by $\{ x_{1},x_{2},\ldots,x_{n}\}$, which is the fundamental group of the $n$-punctured disc $D_{n}$. It is known that the braid group $B_{n}$ is identified with the relative mapping class group $MCG(D_{n},\partial D_{n})$, which is the group of isotopy classes of homeomorphisms of $D_{n}$ which fixes $\partial D_{n}$ pointwise \cite{b}.

The action of the braid groups on $D_{n}$ induces the representation 
$\Phi: B_{n} \rightarrow  \aut(\pi_{1}(D_{n})) = \aut(F_{n})$, which is explicitly written as  
\[ 
\Phi(\sigma_{k})\: :
\left\{
\begin{array}{rll}
x_{k} & \mapsto & x_{k}^{-1}x_{k+1}x_{k} \\
x_{k+1} & \mapsto & x_{k} \\
x_{i} & \mapsto & x_{i} \;\;\;\; (i \neq k,k+1).
\end{array}
\right.
\]

The representation $\Phi$ is identical with the associated group representation $\rho'_{T_{1},q}$.
It is known that both $\rho_{T_{1},q}$ and $\rho'_{T_{1},q}$ are faithful.
\end{exam}

First we show that the representation $\rho_{Q,h}$ is faithful.
\begin{prop}
For a pointed quandle $(Q,h)$, the representation $\rho_{Q,h}: B_{n} \rightarrow \aut(Q^{*n})$ is faithful.
\end{prop}

\begin{proof}
Let us consider the subquandle $T = (\{ h \},*)$, which is isomorphic to the trivial one quandle $T_{1}$, and consider the subquandle $T^{*n} \subset Q^{*n}$. By definition, the restriction $\rho_{Q,h}|_{T^{*n}}$ coincide with the representation $\rho_{T_{1},h}$.
Since $\rho_{T_{1},h}$ is faithful, $\rho_{Q,h}|_{T}$ is also faithful.
\end{proof}

Now using the representation associated to a pointed quandle $(Q,h)$, we construct the quandle invariant functor.
For an $n$-braid $\beta$, we define a {\it quandle-valued invariant} $I_{\beta}(Q,h)$ as the quotient of $Q^{*n}$ by the set of relations $ \{ [\rho_{Q,h}(\beta)] ( q ) = q \: | \: q \in Q^{*n} \}. $
For a pointed quandle morphism $f: (Q,h) \rightarrow (R,i)$, we define a morphism between quandle-valued invariants $I_{\beta}(f): I_{\beta}(Q,h) \rightarrow I_{\beta}(R,i)$ by $[I_{\beta}(f)](q_{i}) = [f(q)]_{i}$. This defines a functor $I_{\beta}:\mathcal{PQ} \rightarrow \mathcal{Q}$. 

\begin{thm}
The functor $I_{\beta}$ is a knot invariant. 
\end{thm}
\begin{proof}
Recall that the Markov theorem (see \cite{b}, for example) states that the closure of two braids $\alpha, \beta$ represent the same oriented link if and only if $\alpha$ is converted to $\beta$ by applying following two operations.
\begin{description}
\item[Conjugation] $\alpha \rightarrow \gamma^{-1}\alpha\gamma$ where $\alpha,\gamma \in B_{n}$.
\item[(De)Stabilization] $\alpha \leftrightarrow \alpha\sigma_{n}^{\pm 1}$ where $\alpha \in B_{n}$.
\end{description}

First we show the invariance under the conjugation. Since $\rho_{Q,h}(\gamma)$ is an automorphism of $Q^{*n}$,
the set of relations $\{ [\rho_{Q,h}(\alpha)](q) = q \}$ is equivalent to the set of relations $\{ [\rho_{Q,h}(\alpha\gamma)](q) = [\rho_{Q,h}(\gamma)](q)\}$.  Hence $I_{\alpha}(Q,h)$ and $I_{\gamma^{-1} \alpha \gamma }(Q,h)$ are isomorphic as a quandle.

The isomorphism $\tau_{\gamma}$ between $I_{\alpha}(Q,h)$ and $I_{\gamma^{-1}\alpha\gamma}(Q,h)$ is given by $\tau_{\gamma}(q_{i}) = [\rho_{Q,h}(\gamma)](q_{i})$. Thus, the following diagram commutes for any pointed quandle morphisms $f: (Q,h) \rightarrow (R,i)$.
\[ \xymatrix{
I_{\alpha}(Q,h)\; \ar[r]^{I_{\alpha}(f)} \ar[d]^{\tau_{\gamma}} \;& I_{\alpha}(R,i) \ar[d]_{\tau_{\gamma}}\\
I_{\gamma^{-1} \alpha \gamma}(Q,h) \;\;\ar[r]^{I_{\gamma^{-1} \alpha \gamma}(f)} \;\;&\; I_{\gamma^{-1} \alpha \gamma }(R,i)
}
\]
Therefore we conclude that $I_{\beta}$ is invariant as a functor under conjugations.

Next we show the invariance under the positive stabilization. First observe that $[\rho_{Q,h}(\alpha)](q_{n+1})=q_{n+1}$.
From the relation $[\rho_{Q,h}(\alpha\sigma_{n})](q_{n+1}) = [\rho_{Q,h}(\alpha)]( q_{n}\os h_{n}) = q_{n+1}$, we obtain the equation
\begin{eqnarray*}
[\rho_{Q,h}(\alpha\sigma_{n})](q_{n}) &  =  & [\rho_{Q,h}(\alpha)] (q_{n+1}* h_{n}) = q_{n+1}* [\rho_{Q,h}(\alpha)](h_{n}) 
\\
& = & [\rho_{Q,h}(\alpha)] (q_{n} \os h_{n}) * [\rho_{Q,h}(\alpha)](h_{n}) \\
& = & [\rho_{Q,h}(\alpha)](q_{n})
\end{eqnarray*}
Thus, the relation $[\rho_{Q,h}(\alpha\sigma_{n})](q_{n}) = q_{n}$ implies $[\rho_{Q,h}(\alpha)] (q_{n})=q_{n}$. Similarly, the relation $[\rho_{Q,h}(\alpha)](q_{n})=q_{n}$ implies $[\rho_{Q,h}(\alpha\sigma_{n})](q_{n}) = q_{n}$. Thus, there is a natural isomorphism $\tau_{+}: I_{\alpha \sigma_{n}}(Q,h) \rightarrow I_{\alpha}(Q,h)$, which is defined by $\tau_{+}(q_{i})= q_{i}$ $(i=1,\ldots,n)$ and $\tau_{+}(q_{n+1}) = q_{n}$. Now the following diagram commutes for any pointed quandle morphism $f:(Q,h) \rightarrow (R,i)$, hence the functor $I_{\beta}$ is invariant under the positive stabilization.
\[ \xymatrix{
I_{\alpha\sigma_{n}}(Q,h)\; \ar[r]^{I_{\alpha\sigma_{n}}(f)} \ar[d]^{\tau_{+}} & I_{\alpha \sigma_{n}}(R,i) \ar[d]_{\tau_{+}}\\
I_{ \alpha}(Q,h) \;\ar[r]^{I_{ \alpha}(f)} & I_{\alpha }(R,i)
}
\]

The invariance under the negative stabilization is similar.
\end{proof}

Now we obtain the first definition of quandle invariant functor.
\begin{defn}
Let $K$ be an oriented knot represented as closed braid $\widehat{\alpha}$. The {\it quandle invariant functor} $I_{K}$ is a functor $I_{\alpha}: \mathcal{PQ} \rightarrow \mathcal{Q}$. For a pointed quandle $(Q,h)$, we call a quandle $I_{K}(Q,h)$ {\it the quandle invariant} associated to $(Q,h)$.
\end{defn}

From definition, it is easy to see the quandle invariant functor $I_{K}$ has the following properties.

\begin{prop}
Let $K$ be an oriented knot.
\begin{enumerate}
\item If $\tau: (Q,h) \rightarrow (R,i)$ is a surjective morphism of pointed quandles, then $I_{K}(\tau): I_{K}(Q,h) \rightarrow I_{K}(R,i)$ is also surjective.
\item For each pointed quandle $(Q,h)$, $I_{K}(Q,h) = I_{-K^{*}}(Q,h)$, where $-K^{*}$ is the dual of $K$.
\end{enumerate}
\end{prop}
\begin{proof}
The assertion (1) is obvious from the definition of $I_{K}(Q,h)$.
Let $K = \widehat{\beta}$. Then, $-K^{*} = \widehat{\beta^{-1}}$, so the relation $[\rho_{Q,h}(\beta)](q)=q$ is equivalent to the relation $[\rho_{Q,h}(\beta^{-1})](q)=q$, hence $I_{K}(Q,h)$ is isomorphic to $I_{-K^{*}}(Q,h)$.
\end{proof}

\section{Diagrammatic description of quandle invariants}

We give an alternative definition of the quandle invariant functor by using a knot diagram. This construction is more combinatorial and useful to study the homology of the quandle invariants.

Let $D$ be an oriented knot diagram, which is a projection of a knot on the plane having transverse double points together with the ``over and under" information. We indicate this information by breaking the under-passing segment. Let $\mA(D) = \{A_{1},A_{2},\ldots,A_{m} \}$ be a set of large arcs, which is an connected component of $D$.
Each large arc $A_{i}$ is decomposed to the subarcs $a_{i,1},a_{i,2},\ldots  ,a_{i,k_{i}}$ by removing the double points of $D$. 
We call these subarcs {\it small arcs} of $D$, and denote the set of small arcs by $\mSA(D)$. For a small arc $a$, we denote the large arc containing $a$ as its subarc by the corresponding large letter $A$.

Let $(Q,h)$ be a pointed quandle and $Q^{*m} = Q_{A}*Q_{B}*\cdots$ be the free product of $m$-copies of $Q$, where $m= \sharp \mA(D)$. Each copy of $Q$ is labeled by the large arc of $D$. For each $q \in Q$ we denote by $q_{A}$ the element of $Q_{A} \subset Q^{*m}$ which corresponds to $q$.

For an element $q \in Q$, we first define the map $c_{q}: \mSA(D) \rightarrow Q^{*m}$ by the following manner. 

\begin{figure}[htbp]
 \begin{center}
\includegraphics[width=45mm]{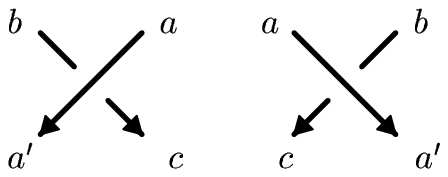}
 \end{center}
 \caption{Labeling of small arcs around the crossing point}
 \label{fig:diagram}
\end{figure}

\begin{enumerate}
\item For a small arc $a$ which contains the starting point of the large arc $A$, we define $c_{q}(a) = q_{A}$.
\item Let $x$ be a crossing point of $D$ and put $a,a',b,c$ as in Figure \ref{fig:diagram}. Assume that we have defined the value $c_{q}(a)$. Then, we define $c_{q}(a')$ by 
\[
\left\{
\begin{array}{l}
c_{q}(a') = c_{q}(a)\os h_{A} \textrm{  if the crossing } $x$ \textrm{ is positive.}\\
c_{q}(a') = c_{q}(a)* h_{A} \textrm{  if the crossing } $x$ \textrm{ is negative.}\\
\end{array}
\right.
\]
\end{enumerate}

Using the map $c_{q}$, we associate a relation $R(x;q)$ for each crossing point $x$ and $q \in Q$ by 
\[ R(x;q) : 
\left\{
\begin{array}{l}
c_{q}(c) = c_{q}(b) * h_{A} \textrm{  if the crossing } $x$ \textrm{ is positive. }\\
c_{q}(c) = c_{q}(b) \os h_{A} \textrm{  if the crossing } $x$ \textrm{ is negative. }\\
\end{array}
\right.
\]

Now we define the quandle invariant $I_{D}(Q,h)$ by the presentation
\[ I_{D}(Q,h) = \langle q_{A} \;\;\; (A \in \mA(D), \; q \in Q) \:|\: R(x;q) \;(x: \textrm{Crossings of }D,\; q \in Q ) \rangle .\] 

As in the definition using braid representation, $I_{D}$ naturally can be seen as a functor $I_{D}:\mathcal{PQ} \rightarrow \mathcal{P}$, by defining $[I_{D}(f)](q_{A}) = f(q)_{A}$ for a pointed quandle morphism $f: (Q,h) \rightarrow (R,i)$.

\begin{prop}
The functor $I_{D}$ is a knot invariant, and it coincides with the quandle invariant $I_{K}(Q,h)$ defined in the previous section.
\end{prop}  

\begin{proof}
We only prove an invariance of the quandle invariant $I_{D}(Q,h)$. Invariance as a functor is routine.
From the definition above, it is easy to confirm that if a diagram $D$ is a closed braid diagram $D = \widehat{\beta}$, then $I_{D}(Q,h)=I_{\beta}(Q,h)$ holds. Thus, only we need to show is that $I_{Q,h}(D)$ is a link invariant. 

\begin{figure}[htbp]
 \begin{center}
\includegraphics[width=80mm]{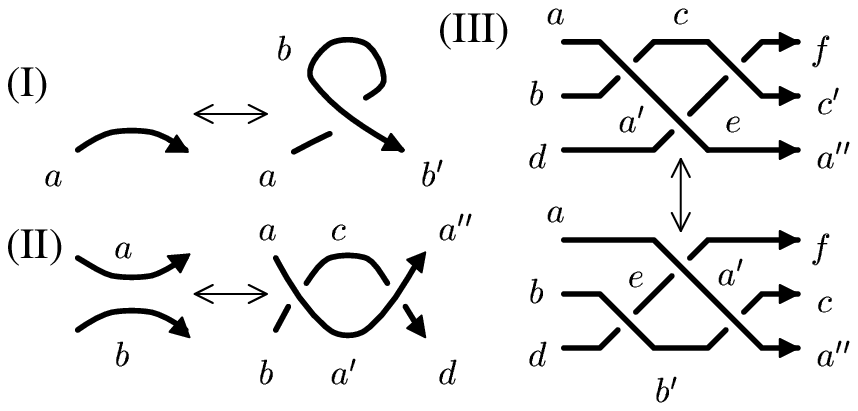}
 \end{center}
 \caption{Reidemeister move invariance}
 \label{fig:rmove}
\end{figure}

{\bf Invariance of Reidemeister move I:}\\

Let $x$ be the newly-added crossing generated by the Reidemeister move I. We consider the case $x$ is a positive crossing, and put $a,b,b'$ as in Figure \ref{fig:rmove}. 
 Then the relation $R(x;q)$ is given by $c_{q}(b)=c_{q}(a) * h_{B}$, hence we obtain $c_{q}(b')=c_{q}(b) \os h_{B} = c_{q}(a)$. Thus this move does not change the quandle invariant. The negative case is similar. \\

{\bf Invariance of Reidemeister move II:}\\

We consider the Reidemeister move II depicted in Figure \ref{fig:rmove} (II). Other cases are similar. Take a small arc $\{a,a',a'',b,c,d\}$ as in Figure \ref{fig:rmove}. 
Then the relations at these two crossings are given by
\[\left\{
\begin{array}{l}
c_{q}(c)= c_{q}(b) * h_{A} \\
c_{q}(d)= c_{q}(c) \os h_{A}.
\end{array}
\right.
\]
Hence we obtain $c_{q}(b)=c_{q}(d)$ and the contributions of the quandle $Q_{C}$ to $I_{D}(Q,h)$ vanish. Hence these two diagrams defines the same quandle. 
\\

{\bf Invariance of Reidemeister move III:}\\

 Take small arcs $\{ a,a',a'',b,b',c,c',d,e,f\}$ as in Figure \ref{fig:rmove} (III).
In the above diagram, three crossings provide the relations 
\[
\left\{
\begin{array}{l}
c_{q}(c)= c_{q}(b)*h_{A} \\
c_{q}(e)= c_{q}(d)*h_{A} \\
c_{q}(f)= c_{q}(e)*h_{C} \\
\end{array}
\right.
\]

Thus, we obtain $c_{q}(f)=(c_{q}(d)*h_{A})*h_{C}$. By putting $q=h$, we obtain $h_{C}=h_{B}*h_{A}$, so 
$c_{q}(f)=(c_{q}(d)*h_{A})*(h_{B}*h_{A}) = (c_{q}(d)*h_{B})*h_{A}$.
Similarly, 
$c_{q}(c')= (c_{q}(b)*h_{A}) \os (h_{B}*h_{A}) = (c_{q}(b) \os h_{B}) * h_{A}$.
On the other hand, from the diagram below, we obtain the relations
\[
\left\{
\begin{array}{l}
c_{q}(e)= c_{q}(d)*h_{B} \\
c_{q}(f)= c_{q}(e)*h_{A} \\
c_{q}(c)= c_{q}(b')*h_{A} = (c_{q}(b)\os h_{B}) *h_{A}\\
\end{array}
\right.
\]
so we also obtain $c_{q}(f)=(c_{q}(d)*h_{B})*h_{A}$.
Thus, the map $c_{q}$ takes the same value on each small arc.
Now it is easy to see that these two diagrams defines the same quandle.
\end{proof}

From this diagrammatic definition, it is quite easy to check our invariant quandle functor is indeed an extension of the classical knot quandle from the diagrammatic definition.

\begin{proof}[Proof of Theorem \ref{thm:main} (1)]
Let us take the trivial $1$-quandle $T_{1}$. Then the relation $R(x)$ at the crossing $x$ is $ q_{C} = q_{B}*q_{A} $, which is a relation appeared in the classical Wirtinger presentation of the knot quandles in \cite{j}.
\end{proof}

Now we show that quandle invariants naturally contain the knot quandle as a subquandle.  
\begin{prop}
\label{prop:object}
There are a natural injection of the knot quandle $\iota: Q_{K} \rightarrow I_{K}(Q,h)$ and a natural surjection to the knot quandle $p: I_{K}(Q,h) \rightarrow Q_{K}$. Moreover, $p\circ \iota = id $.
\end{prop}
\begin{proof}
This proposition follows from Theorem \ref{thm:main} (1) and the fact that the trivial $1$-quandle is the initial and the final object of the category $\mathcal{PQ}$. More precisely, let $i:T_{1} \rightarrow (Q,h)$ be the natural inclusion and $\pi: (Q,h) \rightarrow T_{1}$ be the natural surjection. Then $\iota = i_{*}$ and $p=\pi_{*}$.
\end{proof}

We also remark that from the presentation of $I_{K}(Q,h)$, for each long arc $A$ of a knot diagram $D$, the natural map $Q_{A} \rightarrow I_{D}(Q,h)$ is injection so we may also regard $Q$ as a subquandle of $I_{D}(Q,h)$ by specifying an arc $A$ of the diagram $D$. In the knot theory view point, this corresponds to choose a base point $* \in K$, thus this is equivalent to consider the corresponding long knot.

\section{Homology and cohomology of quandle invariants}

In this section, we study the homology and cohomology groups of the quandle invariant $I_{K}(Q,h)$.
For a quandle $(X,*)$, let $C^{R}_{n}(X)$ be the free abelian group generated by $n$-tuples of elements of $X$ and $C^{D}_{n}(X)$ be the subgroup of $C^{R}_{n}(X)$ generated by $n$-tuples $(x_{1},x_{2},\ldots,x_{n})$ of $X$ with $x_{i}=x_{i+1}$ for some $i$. 
Let $\partial_{n}: C^{R}_{n} \rightarrow C^{R}_{n-1}$ be a homomorphism defined by
\begin{eqnarray*}
\partial_{n}(x_{1},x_{2},\ldots,x_{n}) & = &\sum_{i=1}^{n}(-1)^{i} \left[ (x_{1},x_{2},\ldots, \right. \widehat{x_{i}},\ldots, x_{n}) \\
& & \left. \hspace{1cm} - (x_{1}*x_{i},x_{2}*x_{i},\ldots,x_{i-1}*x_{i},x_{i+1},\ldots, x_{n}) \right]
\end{eqnarray*} 
Then both $(C^{R}_{*}(X),\partial_{*})$ and ($C^{D}_{*}(X),\partial_{*}$) are chain complexes.
Let $C^{Q}_{*}(X)$ be the quotient complex $(C^{R}_{*}(X),\partial_{*})\slash(C^{D}_{*}(X),\partial_{*})$.
For an abelian group $G$,
The $G$-coefficient $n$-th quandle homology and cohomology groups are defined by 
\[ H_{n}^{Q}(X;G) = H_{n}(C_{*}^{Q}(X) \otimes G), \; H^{n}_{Q}(X;G) = H^{n}(\textrm{Hom} (C_{*}^{Q}(X),G))\]
respectively.
The quandle homology is an invariant of quandles, so the homology and cohomology groups of quandle invariant $I_{K}(Q,h)$ also define knot invariants.

\begin{rem}
The above definition of the quandle (co)homology is the simplest one. There are general theory of quandle (co)homologies, including twisted coefficients \cite{ag},\cite{cegs}. Many results in this section also remains true for such generalized homology theories. In particular, we can also extends the ``generalized" cocycle invariants defined in \cite{cegs}, which is an extension of classical cocycle invariants. We mainly restrict the classical (abelian coefficient) case for the sake of simplicity.
\end{rem}

First we observe that the quandle invariant  contains all information of the homology and cohomology of knot quandles.

\begin{lem}
Let $K$ be an oriented knot and $(Q,h)$ be a pointed quandle. Let $\iota :Q_{K} \hookrightarrow I_{K} (Q,h) $ and $p: I_{K}(Q,h) \rightarrow Q_{K}$ be the natural maps in Proposition \ref{prop:object}.
Then for any coefficient group $G$, $\iota$ induces an injection of the homology groups 
$\iota_{*}: H_{*}^{Q}(Q_{L};G) \hookrightarrow H_{*}^{Q}(I_{L}(Q,h);G)$. Similarly, $p$ induces an injection of cohomology groups.
\end{lem} 
\begin{proof}
Since $p\circ \iota = id$, this is clear.
\end{proof}

Now we determine the $1$st quandle homology group $H_{1}^{Q}(I_{K}(Q,h),\Z)$.

\begin{proof}[Proof of Theorem \ref{thm:main} (2)]
Let $D$ be an oriented knot diagram which represents $K$, and take large arcs $A,B \in \mA(D)$.
Since $K$ is a knot, by the definition of quandle invariants there is an element $h_{A,B} \in Q_{K} \subset I_{K}(Q,h) $ such that $q_{A} * h_{A,B} = q_{B}$. Thus $q_{A}$ and $q_{B}$ represents the same $1$st homology class. Hence the map $Q \rightarrow I_{K}(Q,h)$ defined by $q \mapsto q_{A}$ induces an isomorphism of the $1$st quandle homologies.
\end{proof}

Next we proceed to study the $2$nd homology group. 
For an oriented knot $K$,  Eisermann showed that if $K$ is not an unknot, then $H_{2}^{Q}(Q_{K};\Z)=\Z$ and $H_{2}^{Q}(Q_{K};\Z)$ is generated by the {\it orientation class} $[K]$ and if $K$ is unknot, then $H_{2}^{Q}(Q_{K};\Z)=0$ \cite{e}.
 
Now we extend the orientation class for the quandle invariant $I_{K}(Q,h)$.
From now on, we assume that the quandle $Q$ is finite.

\begin{figure}[htbp]
 \begin{center}
\includegraphics[width=45mm]{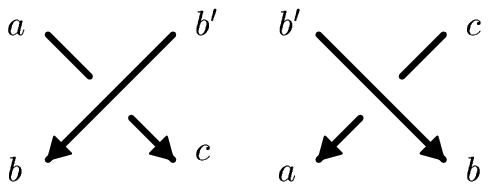}
 \end{center}
 \caption{Definition of 2-chain $(D)$}
 \label{fig:2chain}
\end{figure}

Let $D$ be a knot diagram. Let us define $2$-chain $(D) \in C_{2}^{Q}(I_{K}(Q,h);\Z)$ by   
\[ (D) = \sum_{ q \in Q} \sum_{x} \varepsilon(x) \cdot \{ (c_{q}(a), h_{B}) - (c_{q}(b),h_{B})\}.\]
where $x$ runs all crossing points of $D$ and $\varepsilon(x)$ denotes the sign of the crossing $x$.
For each crossing $x$, we take small arcs $a,b,b',c$ as in Figure \ref{fig:2chain}.

\begin{lem}
\label{lem:orientclass}
The 2-chain $(D)$ is a cycle and its representing homology class $[D] \in H_{2}^{Q}(I_{K}(Q,h);\Z)$ is a knot invariant.
\end{lem}

\begin{proof}
First we show $(D)$ is cycle.
The boundary of $(D)$ is given by
\begin{eqnarray*}
\partial (D) & = & \sum_{x} \sum_{q} \partial \varepsilon(x)\{ (c_{q}(a),h_{B}) - (c_{q}(b),h_{B}) \} \\
 & = & \sum_{x}\sum_{q}\varepsilon(x)(c_{q}(a)*h_{B} - c_{q}(a)) - \varepsilon(x)(c_{q}(b)- c_{q}(b)*h_{B}) \\
 & = & \sum_{x}\sum_{q}\varepsilon(x)( c_{q}(c) -c_{q}(a)) - \sum_{q} \varepsilon(x)(c_{q}(b)-c_{q}(b)*h_{B})\\
 & = & \sum_{x}\sum_{q}\varepsilon(x)(c_{q}(c)-c_{q}(a)) = \sum_{x}\sum_{q}(q_{C}-q_{A}).
\end{eqnarray*}

Hence each long arc $A$ contributes $\partial(D)$ by $\sum_{q}q_{A}$ at its initial point and by $-\sum_{q}q_{A}$ at its end point. Thus, these two contributions cancel each other, so $\partial(D)=0$.

Now we check that its homology class $[D]$ is a knot invariant. 

First we check the invariance of Reidemeister move I. Let $x$ be the newly-added crossing generated by Reidemeister move I. We consider the case $x$ is positive crossing. The negative case is similar. Take a small arc $a,b,b'$ around $x$ as in the Figure \ref{fig:rmove}.
 Then the contribution of the newly-added crossing to the cycle $(D)$ is
\[ \sum_{q} (c_{q}(a),h_{B}) - (c_{q}(b'),h_{B}) = \sum_{q}(c_{q}(a),h_{B}) - (c_{q}(a),h_{B}) = 0,\]
hence the cycle $(D)$ itself is invariant under the Reidemeister move I.

Next We consider the Reidemeister move depicted in Figure \ref{fig:rmove} (II). Other case of Reidemeister move II are similar. Take a small arc $a,a',a'',b,c,d$ as in the figure \ref{fig:rmove}.
Then newly-added two crossings contribute the cycle $(D)$ by
\begin{eqnarray*}
& & \sum_{q}(c_{q}(b),h_{A}) -(c_{q}(a'),h_{A}) - (c_{q}(d),h_{A}) + (c_{q}(a''),h_{A}) \\
& = &  \sum_{q}-(c_{q}(a'),h_{A}) + (c_{q}(a''),h_{A}) \: = \:0 
\end{eqnarray*}
hence the cycle $(D)$ itself is invariant under the Reidemeister move II.

Finally we show the invariance under the Reidemeister move III.
Take small arcs $\{ a,a',a'',b,b',c,c',d,e,f\}$ as in Figure \ref{fig:rmove} (III).
Then the three crossings in the diagram above contribute the cycle by 
\[
\begin{array}{l}
\sum_{q}  \left\{   (c_{q}(b),h_{A}) - (c_{q}(a'),h_{A})  + (c_{q}(d),h_{A}) \right. \\
 \hspace{2cm} \left. - (c_{q}(a''),h_{A}) + (c_{q}(e),h_{C}) -(c_{q}(c'),h_{C})  \right\} 
  \end{array}
\]
 and the three crossings contribute in the diagram below contributes by
\[
\begin{array}{l} 
\sum_{q}  \left\{ (c_{q}(d),h_{B}) - (c_{q}(b'), h_{B}) + (c_{q}(e),h_{A}) \right. \\
\left. \hspace{2cm}- (c_{q}(a'),h_{A}) + (c_{q}(b'),h_{A}) - (c_{q}(a''),h_{A})\right\}.
\end{array}
 \]
Thus their difference is a boundary of the $3$-chain
\[ \sum_{q} \left\{ (q_{D}, h_{B},h_{A}) - (q_{B},h_{B},h_{A}) \right\}.\]
Therefore we conclude that the cycles represent the same homology class.
\end{proof}

The $(Q,h)$-{\it fundamental class} (or orientation class) of a knot $K$ is, by definition, $[(D)] \in H_{2}^{Q}(I_{K}(Q,h);\Z) $ where $D$ is a knot diagram representing $K$. From Lemma \ref{lem:orientclass}, this is independent of the choice of $D$. We denote this homology class by $[K]_{Q,h}$. 
If $Q = T_{1}$, the fundamental class is the same as the orientation class $[K]$ defined by Eisermann \cite{e}.

Now we prove the non-vanishing results and the duality of the $(Q,h)$-fundamental class.

\begin{proof}[Proof of Theorem \ref{thm:main} (3), (4)]
Let $p_{*}: H_{2}^{Q}(I_{K}(Q,h);\Z) \rightarrow H_{2}^{Q}(K_{Q};\Z)$ be the map induced by the natural map $p: I_{Q,h}(K) \rightarrow Q_{K}$ in Proposition \ref{prop:object}. From the definition of $[K]_{Q,h}$, we obtain $p_{*}([K]_{Q,h}) = \sharp\, Q \cdot [K]$.
Since $[K]=0$ if and only if $K$ is an unknot \cite{e}, so we conclude $[K]_{Q,h}=0$ if and only if $K$ is an unknot. 
The assertion (4) follows from the definition of the cycle $(D)$.
\end{proof}

\begin{rem}
The isomorphic class of the quandle invariant $I_{K}(Q,h)$ only depends on the quandle homology class $[h] \in H_{1}^{Q}(Q;\Z)$, because if $h$ and $h'$ represent the same homology class, then the set of the relations coincide. In particular, if $H_{1}^{Q}(Q;\Z)= \Z$, then the isomorphic class of $I_{K}(Q,h)$ does not depends on the choice of the point $h \in Q$.
However, we need to fix a point $h \in Q$ to obtain a functor $I_{K}$ and to obtain the $(Q,h)$-fundamental class $[K]_{Q,h}$.  
\end{rem}

We remark that for unknot $K$, the quandle invariant $I_{K}(Q,h)$ is isomorphic to $Q$, so its $2$nd homology group does not always vanish whereas $H_{2}^{Q}(K_{Q};\Z) = 0$. 

The $(Q,h)$-fundamental class is decomposed as the sum of the {\it partial fundamental class} as follows.
For an element $q \in Q$, we denote by $[q]$ the $h$-orbit of $q$. That is, $[q]$ is a subset of $Q$ defined by $[q]=\{(\cdots (q*h)*h)\cdots *h)*h ,(\cdots (q \os h) \os h)\cdots \os h)\os h \}$.
Then the quandle $Q$ is decomposed as a disjoint union of $h$-orbits as $Q=[h] \coprod [q_{1}] \coprod \cdots \coprod [q_{k}]$.
For a knot diagram $D$ and $q \in Q$, define 
\[ (D)_{[q]} = \sum_{ q \in [q]} \sum_{x} \varepsilon ( x )\cdot (c_{q}(a_{x}),c_{q}(b_{x}) ).\]
where $x$ runs all crossings of $D$. 
As in the proof of Lemma \ref{lem:orientclass}, $(D)_{[q]}$ is also a cycle and its homology class is also independent of a choice of the knot diagram. Thus this also defines a knot invariant $[K]_{Q,h;[q]}$. We call this homology class the {\it partial } $(Q,h)$- {\it fundamental class} of $K$ relative to $[q]$.

From the definition, we obtain the following results, which extends the results for $(Q,h)$-fundamental class.

\begin{cor}
Let $K$ be an oriented knot and $(Q,h)$ be a finite pointed quandle. 
We denote the $h$-orbit decomposition of $Q$ by $Q= [q_{0}] \coprod \cdots \coprod [q_{k}]$. Then,
\begin{enumerate}
\item $[K]_{Q,h} = \sum_{i=0}^{k} [K]_{Q,h;[q_{i}]}$.
\item For each $q \in Q$, $[K]_{Q,h;[q]}$ is trivial if and only if $K$ is unknot.
\end{enumerate}
\end{cor}

We close this section by giving a question about (co)homology group of quandle invariants.
As we have seen, the $1$st homology group of quandle invariants contains no information of $K$, because it is determined by only $Q$. 
 We would like to ask this phenomenon always occurs for all degrees.
\begin{ques*}\textrm{ }
\begin{enumerate}
\item The (co)homology group of $I_{K}(Q,h)$ are always determined by $H_{*}^{Q}(Q_{K};\Z)$ and $H_{*}^{Q}(Q;\Z)$ ?
\item The betti numbers of $I_{K}(Q,h)$ are always determined by the betti number of $Q_{K}$ and $Q$ ?
\end{enumerate}
\end{ques*}

\section{Quandle cocycle invariant via quandle invariant $I_{K}(Q,h)$ }

In this section we extend the quandle cocycle invariants using the quandle invariant $I_{K}(Q,h)$ and study its properties. As in the previous section, we always assume that every pointed quandle $(Q,h)$ is finite.

Let $D$ be a link diagram, and $X$ be a finite quandle. Take a $G$-coefficient quandle  2-cocycle of $X$, $\phi: X\times X \rightarrow G$.
We call a quandle morphism $\rho: I_{K}(Q,h) \rightarrow X$ {\it a }$(Q,h)$-{\it extended} $X$-{\it coloring}. For the classical knot quandle, we call a quandle morphism $\rho: Q_{K} \rightarrow X$ a $X$-{\it coloring}. For each $q \in Q$ and a $(Q,h)$-extended coloring $\rho$, we denote by $\rho_{q}$ the map $\rho\circ \overline {c_{q}}: \mSA(D) \rightarrow X$, where $\overline{c_{q}}$ is the composite of the coloring map $c_{q}: \mSA(D) \rightarrow Q^{* \sharp m}$ and the projection map $p: Q^{* m} \rightarrow I_{K}(Q,h)$.

For each crossing $x$ of $D$, we define a weight $W(x,q;\rho)$ by  
\[
W(x,q;\rho) = \{ \phi(\rho_{q}(a), \rho(h_{B}) ) \cdot \phi(\rho_{q}(b),\rho(h_{B}))^{-1}\}^{\varepsilon(x)}. 
\] 

The $(Q,h)$-{\it extended quandle cocycle invariant} is defined by the sum of the all weights 
\[ \Phi_{\phi,(Q,h)} (D) = \sum_{\rho} \prod_{q\in Q}\prod_{x} W(x,q;\rho) \in \Z[G]\]
where $x$ runs all crossings of $D$ and $\rho$ runs all $Q$-extended $X$-colorings. 

\begin{thm}
The $(Q,h)$-extended quandle cocycle invariant is equal to the value 
$\sum_{\rho} \langle [\phi],\rho_{*}([D]) \rangle $. Here $\langle\;,\: \rangle$ represents the pairing of the homology and the cohomology.
Thus, 
$\Phi_{\phi,(Q,h)}(D)$ is a knot invariant and its value depends on the cohomology class $[\phi] \in H_{Q}^{2}(X;G)$.
\end{thm}
\begin{proof}
From the definition of $[D]$, we obtain the equality
\[ \langle [\phi] ,\rho_{*}([D])\rangle = \sum_{\rho} \left (\prod_{q \in Q} \prod_{x} \{  \phi(\rho(c_{q}(a)), \rho(h_{B}) ) \cdot \phi ( \rho (c_{q}(b)),\rho (h_{B}) )^{-1}\}^{\varepsilon(x)} \right).  \]
The right hand is the definition of $\Phi_{\phi,(Q,h)}(D)$.
\end{proof}
 
By definition, for the trivial $1$-quandle $T_{1}$, the $T_{1}$-extended quandle cocycle  invariant coincide with the classical quandle cocycle invariant $\Phi_{\phi}(K)$ defined in \cite{cjkls}.

From the homological viewpoint of the cocycle invariant, we can decompose the $(Q,h)$-extended quandle cocycle invariants by using partial fundamental classes.
For pointed quandle $(Q,h)$ and an element $q \in Q$, let us define the {\it partial quandle cocycle invariant} $\Phi_{\phi,(Q,h);[q]}(K)$ by
\[ \Phi_{\phi,(Q,h);[q]}(K) = \sum_{\rho} \langle [\phi],\rho_{*}([K]_{Q,h;[q]}) \rangle. \]

\begin{cor}
Partial quandle cocycle invariant $\Phi_{\phi,(Q,h);q}(K)$ is a knot invariant.
\end{cor}

Now let us proceed to study the properties of $(Q,h)$-extended quandle cocycle invariants. 
Unfortunately, in many cases $(Q,h)$-extended quandle cocycle invariants are determined by $Q$ and usual quandle cocycle invariants $\Phi_{\phi}(K)$ as we shall explain below.

Before stating our result, we review some notions about quandle morphisms. 
We say a quandle morphism $f: Q \rightarrow R$ is {\it trivial} if $f(Q) = \{r\}$ for some $r \in R$. 
For each element $q \in Q$, the map $[*q]: Q \rightarrow Q$, $x \mapsto x*q$ defines a quandle automorphism of $Q$. An inner automorphism group $\Inn (Q)$ of $Q$ is a subgroup of $\aut(Q)$ generated by  $\{[*q] \: | \: q \in Q\}$. By the definition of inner automorphisms, each nontrivial inner automorphism has at least one fixed point.

\begin{figure}[htbp]
 \begin{center}
\includegraphics[width=120mm]{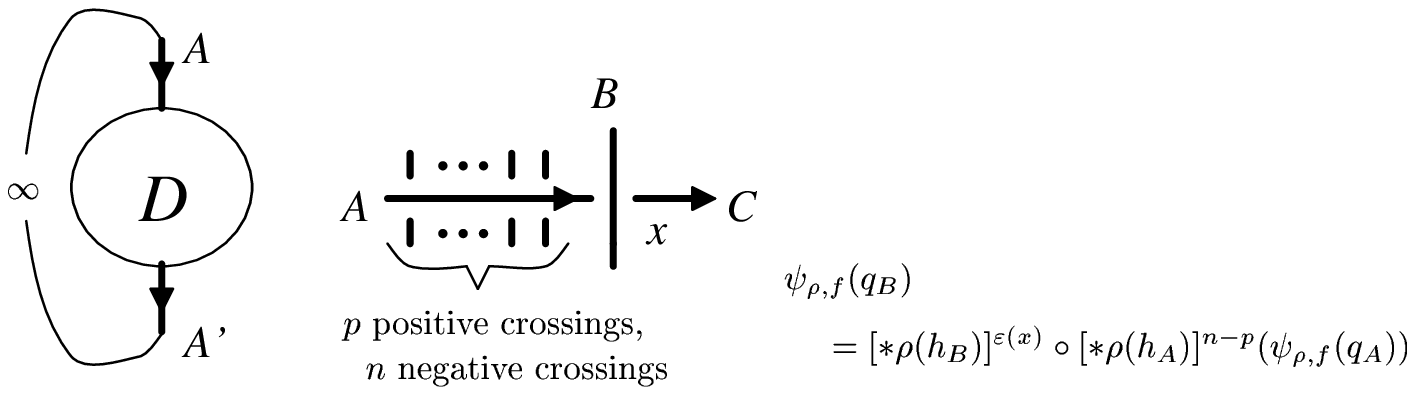}
 \end{center}
 \caption{The definition of $\psi_{\rho,f}(q_{C})$ and $A_{\rho,D}$}
 \label{fig:determinecolor}
\end{figure}

First we observe the relationships between $(Q,h)$-extended quandle coloring $\psi: I_{K}(Q,h) \rightarrow X$ and usual knot coloring $\rho: Q_{K} \rightarrow X$.
Let $D$ be an oriented knot diagram of $K$. We may assume that $D$ is a long knot diagram with special arcs $A$ and $A'$ which contains the point of infinity, as shown in Figure \ref{fig:determinecolor}. Then, $Q_{A}$ is regarded as a subquandle of $I_{K}(Q,h)$.

For a $(Q,h)$-extended coloring $\psi: I_{K}(Q,h) \rightarrow X$, one can obtain quandle morphisms $\rho_{\psi}: Q_{K} \rightarrow X$ and $f_{\psi}: Q \rightarrow X$ which satisfy $\rho_{\psi}(h_{A}) = f_{\psi}(h)$ by considering the restriction of $\psi$ to $Q_{K}$ and $Q_{A}$ respectively.

Conversely, let $\rho: Q_{K} \rightarrow X$ and $f: Q \rightarrow X$ be quandle morphisms which satisfy $\rho(h_{A}) = f(h)$. Now we construct a $(Q,h)$-extended coloring by extending $\rho$ and $f$.
First we define $\psi_{\rho,f}(q_{A}) = f(q)$. Then, using the defining relations of $I_{K}(Q,h)$, we can uniquely determine the value $\psi_{\rho,f}(q_{B})$ for other arcs $B \in \mA(D)$ as the following way.
Let $A,B,C$ be arcs of the diagram $D$ around a crossing point $x$, as in Figure \ref{fig:determinecolor}. 
Assume that we have already defined the value $\psi_{\rho,f}(q_{A})$.
let $p,n$ be the number of positive and negative crossing points contained in $A$.
Then, we define the value $\psi_{\rho,f}(q_{C})$ by 
\[ \psi_{\rho,f}(q_{C}) = [* \rho(h_{B})]^{\varepsilon(x) }\circ [*\rho(h_{A})]^{n-p}(\psi_{\rho,f}(q_{A})). \]

This procedure defines a inner automorphism $A_{\rho,D}: X \rightarrow X$, which sends the color on $A$ to the color on $A'$.
By definition, the inner automorphism $A_{\rho,D}$ only depends on the diagram $D$ and $\rho$.

Now this construction of $\psi_{\rho,f}$ defines a well-defined morphism $I_{K}(Q,h) \rightarrow X$ if and only if the colorings on $A$ and $A'$ coincide, that is, $ A_{\rho,D}(f(q)) = f(q)$ holds for all $q \in Q$.

Summarizing, we proved the following lemma.
\begin{lem}
\label{lem:color} 
There is an one-to-one correspondence between the $(Q,h)$-extended $X$-colorings of knot $K$ and the pair $(\rho,f)$ consisting of a $X$-coloring $\rho:Q_{K} \rightarrow X$ and  a quandle morphism $f: Q \rightarrow X$ which satisfy the two conditions
\begin{enumerate}
\item $\rho(h_{A})=f(h)$.
\item $A_{\rho,D}(f(q)) = f(q)$ for all $q \in Q$.
\end{enumerate}
\end{lem}
 
Now we provide some computations of $(Q,h)$-extended cocycle invariants.

\begin{prop}
\label{prop:cocycleinvariant1}
Let $(Q,h)$ and $X$ be finite (pointed) quandles and 
$K$ be an oriented knot which is non-trivial. Assume that one of the following conditions holds.
\begin{enumerate}
\item There are no non-trivial quandle morphisms $Q \rightarrow X$.
\item Each non-trivial inner automorphism of $X$ has only one fixed point, and $A_{\rho,D}$ is non-trivial for all quandle morphisms $\rho: Q_{K} \rightarrow X$.
\end{enumerate}
Then, 
\[ \Phi_{\phi,(Q,h)} (K) =P^{\sharp Q}(\Phi_{\phi}(K)) \]
holds where $P^{i}: \Z G \rightarrow \Z G$ is a map defined by $g \mapsto g^{i}$ and $\Phi_{\phi}(K)$ is the classical quandle cocycle invariant.
\end{prop}
\begin{proof}

First assume that (1) holds. Since there is no non-trivial quandle morphisms from $Q$ to $X$, by Lemma \ref{lem:color}
there is an one-to-one correspondence between the set of $(Q,h)$-extended colorings and the usual knot colorings. For a usual knot coloring $\rho$, we denote by $\psi_{\rho}$ its corresponding $(Q,h)$-extended coloring. This coloring $\psi_{\rho}$ is defined by $\psi_{\rho}(q_{A}) = \rho(h_{A})$.
Let us denote by $W'(x,\rho)$ the classical weight $\phi(\rho(h_{A}),\rho(h_{B}))^{\varepsilon (x)}$.
Then, the (classical) quandle cocycle invariant of $K$ is defined by the sum of classical weights $\Phi_{\phi}(K) = \sum_{\rho}\prod_{x} W'(x;\rho)$.

From the definition of $\psi_{\rho}$, for a crossing point $x$ of $D$ and an arbitrary element $q \in Q$, the equality $W(x,q;\psi_{\rho}) = W'(x;\rho)$ holds.

Similarly, if (2) holds, by Lemma \ref{lem:color}, there is also an one-to-one correspondence between the set of $(Q,h)$-extended colorings and the usual knot colorings. Thus in this case the same equality of weights holds.

Thus in either case, we obtain
\[ \Phi_{\phi,(Q,h)} (K) = \sum_{\rho}\prod_{q}\prod_{x} W(x,q;\psi_{\rho})  = \sum_{\rho}\prod_{x} W'(x;\rho)^{\sharp Q} =  P^{\sharp Q}(\Phi_{\phi}(K))\]
holds.
\end{proof}

\begin{exam}
We give some examples the criterion of Proposition \ref{prop:cocycleinvariant1} works.  
\begin{enumerate}
\item Let $S_{4}$ be the Alexander quandle $\Z_{2}[T,T^{-1}] \slash T^{2}+T+1$ and $R_{3}$ be the dihedral quandle of order $3$. Then, there are no non-trivial quandle morphism $R_{3} \rightarrow S_{4}$. Hence for all $h \in R_{3}$, the $(R_{3},h)$-extended cocycle invariant for a cocycle of $S_{4}$ is determined by the usual cocycle invariant.
\item Let $K$ be a knot with less than 7 crossings. Then, the inner automorphisms $A_{D,\rho}$ are non-trivial for all coloring map $\rho: Q_{K} \rightarrow S_{4}$. Since each non-trivial inner automorphism of $S_{4}$ has exactly one fixed point, we conclude that for such knots, the $(S_{4},h)$-extended cocycle invariant for a cocycle of $S_{4}$ is determined by the usual cocycle invariant.
\end{enumerate}
\end{exam}

Next we consider the case that $Q$ is a trivial quandle. In this case we can also represent the extended cocycle invariants by the classical cocycle invariants, but the formula is slightly complicated.

\begin{prop}
Let $T_{m}$ be the trivial $m$-quandle and $\phi$ be a $2$-cocycle of a finite quandle $X$.  
Then there exist integers $\{N_{i}\}_{i=1,\ldots,m}$ such that for every oriented knot $K$, $(T_{m},h)$-extended quandle cocycle invariants satisfy the equality
\[ \Phi_{\phi, ( T_{m},h)}(K) = \sum_{i=1}^{m} N_{i}\cdot P^{i}(\Phi_{\phi}(K)).\]
The integers $N_{i}$ depends on only $m$ and $X$.
\end{prop}
\begin{proof}
We always assume $m \geq 2$, because $m=1$ case is trivial.
Let $D$ be a long knot diagram which represents $K$, and fix an arc $A$ so that we may regard $T_{m}=Q_{A} $ as a subquandle of $I_{K}(T_{m},h)$. 
For a quandle morphism $\rho: Q_{K} \rightarrow X$, let $F_{\rho}$ be the set of quandle morphisms from $Q$ to $X$ which satisfy $\rho(h_{A}) = f(h)$ and $X'$ be the image of $\rho(Q_{K})$.

Since the knot quandle $Q_{K}$ is connected (that is, the action of the inner automorphism group of $Q_{K}$ is transitive), $X'$ is also connected. Thus for all $f \in F_{\rho}$ and $q \in Q$, if $f(q) \not \in X'$, then $f(q) * x' = f(q)$ holds for all $x' \in X'$. Similarly, by the same reason, if $f(q) \in X'$, then $f(q) = f(h)$ holds.

By definition, $A_{\rho,D}$ belongs to the subgroup of $\Inn(X)$ generated by $\{[* x'] \: | \: x' \in X'\}$.
Since $\rho$ is a knot coloring, $A_{\rho,D}(\rho(h_{A})) = \rho(h_{A})$ always holds.
Therefore, from the above observations, the inner automorphism $A_{\rho,D}$ is always trivial when it is restricted to $f(T_{m})$. Thus by Lemma \ref{lem:color}, for all $f \in F_{\rho}$, a pair $(\rho,f)$ always defines a $(T_{m},h)$-extended coloring $\psi_{\rho,f}$.

For $q \in Q$, if $f(q) \neq f(h)$, then $\psi_{\rho,f}(c_{q}(b)) = f(q)$ holds for all small arc $b$.
Thus in this case, the weight $W(x,q;\psi_{\rho,f}) = 1$ holds for all crossing $x$ of $D$.
On the other hand, if $f(q) = f(h)$, then $\psi_{\rho,f}(c_{q}(b)) = \rho(h_{B})$ holds for all small arc $b$. Thus, in this case $W(x,q;\psi_{\rho,f})$ is equal to the classical weight $W'(x;\rho)$.

Now we define integer $N_{i}$ as follows. Let us take an element $x \in X$ and define $N_{i}= \sharp \, \{ f: T_{m} \rightarrow X \: | \: \sharp \, f^{-1}(x)= i\}$ for $i=1,\ldots,m$.
The numbers $N_{i}$ are independent of the choice of $x$, and only depend on $m$ and $X$.

Then, by using the obtained equality of weights, we conclude that the equality 
\begin{eqnarray*}
\Phi_{\phi,(T_{m},h)}(D) & = & \sum_{\rho} \sum_{F_{\rho}} \prod_{q} \prod_{x} W(x,q; \psi_{\rho,f})   =  \sum_{\rho} \sum_{i=1}^{m} N_{i} \prod_{x} W'(x;\rho)^{i} \\
 & = & \sum_{i=1}^{m} N_{i} \left( \sum_{\rho}\prod_{x} W'(x;\rho)^{i} \right )  =  \sum_{i=1}^{m} N_{i}\cdot P^{i}(\Phi_{\phi}(D))
\end{eqnarray*}
 holds.
\end{proof}

As these examples suggest, in many simple cases extended quandle cocycle invariants are determined by the usual quandle cocycle invariants (and the pointed quandle $(Q,h)$). In fact, the author cannot find an example of knots whose extended cocycle invariants can distinguish while the corresponding classical cocycle invariant cannot.
Thus, we would like to pose the following question.

\begin{ques*}
Let $X$ be a finite quandle and $\phi$ be a $2$-cocycle of $X$.
For two oriented knots $K$ and $K'$, if their cocycle invariants $\Phi_{\phi}(K)$ and $\Phi_{\phi}(K')$ are the same,
then $(Q,h)$-extended cocycle invariants $\Phi_{\phi,(Q,h)}(K)$ and $\Phi_{\phi,(Q,h)}(K')$ are always the same for all finite pointed quandle $(Q,h)$ ?
\end{ques*}

Even if the above question has a positive answer, it might be difficult to construct an explicit formula to write $(Q,h)$-extended cocycle invariants by the classical cocycle invariants.

\section{Spatial realization of quandle invariants}

In this section we describe a spacial realization of the quandle invariant $I_{K}(Q,h)$ in special case. The contents of this section is a direct extension of section 3 of \cite{cp} and proofs are almost the same, so we only sketch the proof. We remark that this approach does not produce the quandle invariant functor, because it is not known that every quandle and quandle morphisms admits a spatial realizations as  the fundamental quandles.

First we review the notion of the fundamental quandle, introduced by Joyce \cite{j}.
A pointed pair of topological space is a triple $(X,A,*)$ consisting of a topological space $X$, its subspace $A$ and a base point $* \in X \backslash A$. 
A map of pair of topological spaces is a continuous map $f: (X,A,*) \rightarrow (Y,B,*)$ with $f^{-1}(B)=A$ and $f(*)=*$.

\begin{figure}[htbp]
 \begin{center}
\includegraphics[width=60mm]{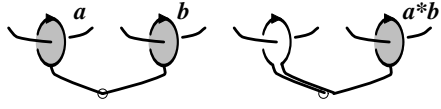}
 \end{center}
 \caption{Quandle operation $*$}
 \label{fig:qoperation}
\end{figure}

Let $N =\{z \in \mathbb{C}\; | \; |z| \leq 1\} \cup \{z \in \mathbb{R} \subset \mathbb{C} \;| \; -5 \leq z \leq -1 \}$. We denote by $\hoop$ the pointed pair of topological space $(N,0,-5)$.
The fundamental quandle $Q(X,A,*)$ of a pointed pair of topological space $(X,A,*)$ is defined as the homotopy classes of the map $f: \hoop \rightarrow (X,A,*)$. The quandle operation $*$ is defined by the Figure \ref{fig:qoperation}.

Under some conditions, for example, in the case $A$ is a codimension two embedding, we can define the notion of positive intersections.
The {\it positive fundamental quandle} $Q^{+}(X,A,*)$ is a subquandle of $Q(X,A,*)$ generated by the map $f:\hoop \rightarrow (X,A,*)$ which positively intersects with $K$ at the point $f(0)$. 

Geometrically, the knot quandle $Q_{K}$ of a knot $K$ is defined as the positive fundamental quandle $Q^{+}(S^{3},K,*)$.

Let $K=\widetilde{\beta}$ be a link represented as a closed $n$-braid.
Let $Q$ be a positive fundamental quandle of a pointed topological pair $(X,A,*)$ and $f: \hoop \rightarrow (X,A,*)$ be a map which represents the element $h \in Q$.

\begin{figure}[htbp]
 \begin{center}
\includegraphics[width=65mm]{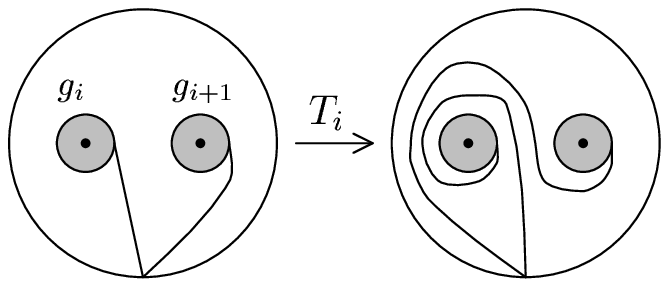}
 \end{center}
 \caption{Maps $g_{i}$ and $T_{i}$}
 \label{fig:map}
\end{figure}

Let $D$ be a 2-disc $D = \{z \in \mathbb{C} \;| \; |z| \leq n+1\}$, $P = \{ p_{i} = (i,0) \in D \;| (i=1,2,\cdots,n) \}$ and $*$ be the base point lying on $\partial D$. We consider the pointed topological pair $(D,P,*)$.
Let $g_{i}: \hoop \rightarrow (D,P,*)$ be the map defined as in Figure \ref{fig:map} and we denote its image in $D$ by $N_{i}$. Now glue $n$-copies of a pointed topological pair $(X,A,*)$ to $(D,P,*)$ along $N_{i}$ by the map $f\circ g_{i}^{-1}$. Let us denote the obtained pointed topological pair by  $(Z,S,*)$.
Then, the fundamental quandle $Q^{+}(Z,S,*)$ is isomorphic to $Q^{*n}$.

Let $C_{i}$ (resp. $C_{i,i+1}$) be a simple closed curve in $D$ which encloses $p_{i}$ (resp. $p_{i}$ and $p_{i+1}$). We denote the half-Dehn twist along $C_{i}$ (resp. $C_{i,i+1}$) by $\tau_{i}$ (resp. $\tau_{i,i+1}$). Let $T_{i}:D \backslash P \rightarrow D \backslash P$ be a homomorphism defined by $\tau_{i}^{-3}\tau_{i+1}^{-1}\tau_{i,i+1}$ (See Figure \ref{fig:map}).
The homeomorphism $T_{i}$ can be extended as a homeomorphism of the pointed topological pair $T_{i}^{X}: (Z,S,*) \rightarrow (Z,S,*)$. 

The following lemma is proved by the same way as in the proof of Proposition 3.2 in \cite{cp}. 
\begin{lem}
\label{lem:geomrep}
Let $(Q,h)$ be a based quandle and $(X,K,*)$ be a pointed topological pair defined as the above. Then the induced homomorphism
$\Phi: B_{n} \rightarrow Aut(Q^{*n})$ defined by $\sigma_{i} \mapsto (T^{X}_{i})_{*}$ is identical with the associated braid representation $\rho_{Q,h}$.
\end{lem}

Now the geometric construction goes as follows. Fix a word representative of $\beta$, and let $B :(Z,S,*) \rightarrow (Z,S,*)$ be a homomorphism which corresponds to $\beta$. 
Then by Lemma \ref{lem:geomrep}, the induced map $B_{*} : Q^{*n} \rightarrow Q^{*n}$ is identical with the image of the associated braid representation $\rho_{Q,h}(\beta)$.
Let $(M(Z),M(S),*)$ be the mapping torus of $B$. Then the total space $M(Z)$ has a torus boundary $\partial D \times S^{1}$.
Along this torus boundary, attach a solid torus so that $\{*\} \times S^{1}$ is identified with $\partial D^{2} \times \{\textrm{point}\}$. We denote the obtained pointed pair of space by
$(\Omega, M(S),*)$.
\begin{thm}
\label{thm:space}
The positive fundamental quandle of the pointed pair of a topological space $(\Omega, M(S), *)$ is isomorphic to the quandle invariant $I_{K}(Q,h)$.
\end{thm} 

\begin{figure}[htbp]
 \begin{center}
\includegraphics[width=50mm]{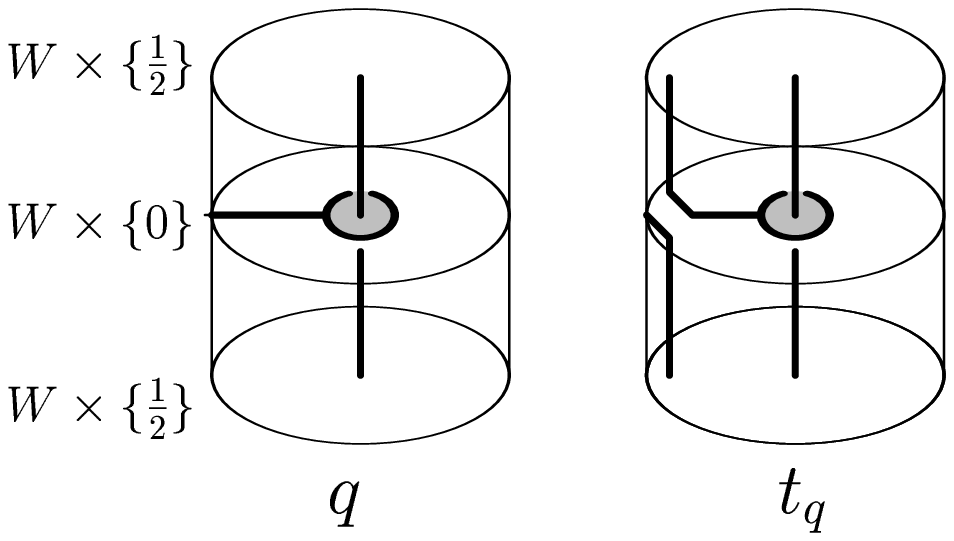}
 \end{center}
 \caption{Generators $q$ and $t_{q}$ of $Q^{+}(M(Z),M(S),*)$ }
 \label{fig:geomgen}
\end{figure}

\begin{proof}
Let $T_{Q}$ be another copy of $Q^{n}$. We denote an element of $T_{Q}$ corresponds to $q \in Q^{*n}$ by $t_{q}$. Then the positive fundamental quandle $Q^{+}(M(Z),M(S),*)$ has a presentation 
\[ Q^{+}(M(Z),M(S),*) = \langle q , t_{q} \: | \: t_{q} =  [\rho_{Q,h}(\beta)](q) \rangle \]
Geometrically, $t_{q}$ is represented by a map depicted in Figure \ref{fig:geomgen}.
Then, gluing a solid torus corresponds to the adding relations $\{q = t_{q}\}$, hence the positive fundamental quandle of $(\Omega, M(S), *)$ has a presentation
\[ Q^{+}(\Omega,M(S),*) = \langle q \in Q^{*n} \: | \:  q =  [\rho_{Q,h}(\beta)](q) \rangle,\]
 which is a presentation of the quandle invariant $I_{K}(Q,h)$.
\end{proof}

This points of view provides a geometrical meaning of quandle invariants in some special cases.
\begin{exam}
Let $FQ_{n}$ be the rank $n$ free quandle generated by $\{q_{1},\ldots, q_{n}\}$, which is the positive fundamental quandle of $(D^{2},\{p_{1},\ldots, p_{n} \},*)$.
By Theorem \ref{thm:space}, $I_{K}(FQ_{n},q_{i})$ is isomorphic to the link quandle of the $n$-parallelization of the knot $K$.
\end{exam}

\end{document}